\documentclass[review]{elsarticle}

\usepackage{lineno,hyperref}
\modulolinenumbers[5]

\journal{Indagationes Mathematicae}

\usepackage{amsfonts,amsmath,latexsym,amssymb}
\newtheorem{theorem}{Theorem}[section] 
\newtheorem{lemma}[theorem]{Lemma}     

\newtheorem{remark}[theorem]{Remark}

\begin{document}

\begin{frontmatter}

\title{Probability Theoretic Generalizations of Hardy's and Copson's Inequality}

\author{Chris A.J. Klaassen \\
   Korteweg-de Vries Institute for Mathematics \\
   University of Amsterdam\\
   Netherlands\\
   c.a.j.klaassen@uva.nl \\
\vspace{1cm}   
  {\em Dedicated to centenarian Jaap Korevaar} }

\begin{abstract}
A short proof of the classic Hardy inequality is presented for $p$-norms with $p>1$.
Along the lines of this proof a sharpened version is proved of a recent generalization of Hardy's inequality in the terminology of probability theory.
A probability theoretic version of Copson's inequality is discussed as well.
Also for $0<p<1$ probability theoretic generalizations of the Hardy and the Copson inequality are proved.
\end{abstract}


\end{frontmatter}

\linenumbers

\section{The classic Hardy inequality for $p>1$}
As described in detail in \cite{MR2256532} the classic Hardy inequality was developed in the years 1906--1928.
Its integral version may be formulated as follows.
If $p>1$ and $\psi$ is a nonnegative measurable function on $(0,\infty)$, then
\begin{eqnarray}\label{continuousHardy}
\left\{\int_0^\infty \left ( \frac{1}{x} \int_0^x \psi (y) dy \right )^p \, dx \right\}^{1/p}
\leq \frac p{p-1} \left\{ \int_0^{\infty} \psi^p(y)\, dy \right\}^{1/p}
\end{eqnarray}
holds.
By Lebesgue's differentiation theorem, Tonelli's theorem and H\"older's inequality we have
\begin{eqnarray}\label{proofcontinuousHardy}
\lefteqn{ \int_0^\infty \left ( \frac{1}{x} \int_0^x \psi \right )^p \, dx
= \int_0^\infty x^{-p} \int_0^x p \left(\int_0^y \psi \right)^{p-1} \psi(y)\, dy \, dx }\\
&& = p \int_0^\infty \int_y^\infty x^{-p} \, dx \left( \int_0^y \psi \right)^{p-1} \psi(y)\, dy
 = \frac p{p-1} \int_0^\infty \left( \frac 1y \int_0^y \psi \right)^{p-1} \psi(y) \, dy  \nonumber \\
&& \leq \frac p{p-1} \left[\int_0^\infty \left( \frac 1y \int_0^y \psi \right)^p dy \right]^{(p-1)/p}
\left[ \int_0^\infty \psi^p \right]^{1/p}, \nonumber
\end{eqnarray}
which implies (\ref{continuousHardy}).
This proof is a smoothed version of the proof of \cite{Hardy:1925}, whose application of partial integration in stead of Tonelli's theorem introduced some technical complications; see Section 8 of \cite{MR2256532} and also observe that partial integration can always be done by applying Tonelli's or Fubini's theorem.

With $\psi(y) = \sum_{k=1}^\infty c_k {\bf 1}_{[k-1,k)}(y)$ and $c_1 \geq c_2 \geq ... \geq 0$ we have
\begin{eqnarray*}
\frac 1x \int_0^x \psi(y)\,dy = \frac{\sum_{k=1}^{n-1} c_k + (x-n+1)c_n}x \geq \frac 1n \sum_{k=1}^n c_k,\ x \in [n-1,n],
\end{eqnarray*}
and by substituting this into (\ref{continuousHardy}) (cf. \cite{Hardy:1925} and Section 8 of \cite{MR2256532}) we obtain
 the sequence version of Hardy's inequality, which states that for $p>1$ and nonnegative $c_1, c_2, ...$
\begin{equation}\label{discreteHardy}
\sum_{n=1}^\infty \left ( \frac{1}{n} \sum_{k=1}^n c_k \right )^p \leq \left ( \frac{p}{p-1} \right )^p  \sum_{k=1}^\infty c_k^p
\end{equation}
holds.
Note that rearranging the $c_k$ in nonincreasing order does not change the value of $\sum_{k=1}^\infty c_k^p$,
but it might increase the value of $1/n \sum_{k=1}^n c_k$ for some values of $n$ without decreasing this value for any $n$,
as may be seen by interchanging any two $c_k$ that are not in decreasing order.

\section{A probability theoretic Hardy inequality for $p>1$}
Both the integral and sequence version of the Hardy inequality can be derived from a probability theoretic generalization presented in \cite{MR4346672}, which we sharpen as follows.
\begin{theorem}\label{Hardytheorem}
Let $X$ and $Y$ be independent random variables with distribution function $F$ on $({\mathbb R}, {\cal B})$,
and let $\psi $ be a measurable function on $({\mathbb R}, {\cal B})$, not identically 0.
Furthermore, let $p>1$ and let $\alpha \in [0,1]$ be the unique root of
\begin{equation}\label{root}
E\left( |\psi(Y)| \right) \left[ p - 1 + \alpha \right] - p \left\{ E\left( |\psi(Y)|^p \right) \right\}^{1/p} \alpha^{1/p} = 0.
\end{equation}
Then
\begin{equation}\label{probabilityHardy}
\left\{ E\left( \left| \frac{E\left( \psi(Y) {\bf 1}_{[Y \leq X]} \mid X \right)}{F(X)} \right|^p \right) \right\}^{1/p}
\leq \frac p{p - 1 + \alpha} \left\{ E\left( |\psi(Y)|^p \right) \right\}^{1/p}
\end{equation}
holds.
\end{theorem}

In (\ref{probabilityHardy}) $E( \cdot \mid X)$ denotes conditional expectation given $X$.
Note that this inequality (\ref{probabilityHardy}) is trivial if $\psi$ is not $p$-integrable.
Furthermore, (\ref{probabilityHardy}) with $\psi$ nonnegative implies (\ref{probabilityHardy}) for arbitrary $\psi$, since the absolute value of the (conditional) expectation of a random variable is bounded from above by the (conditional) expectation of the absolute value of this random variable.

If $\psi(Y)$ equals a constant a.s., then we have $\alpha=1$ and (\ref{probabilityHardy}) is an equality.
For $p=2$ the root from (\ref{root}) equals $\alpha = (\sqrt{E\psi^2(Y)} - \sqrt{{\rm var~}|\psi(Y)|})^2 (E|\psi(Y)|)^{-2}$.

Note that taking $F$ uniform on $(0,K)$, denoting the root of (\ref{root}) by $\alpha_K$,
multiplying (\ref{root}) by $K^{1/p}$ and (\ref{probabilityHardy}) by $K,$
and taking limits as $K \to \infty$ we see that $\alpha_K$ converges to 0 and we arrive at Hardy's inequality (\ref{continuousHardy}).
Analogously taking $F$ uniform on $\{1,  ..., K\}$ we obtain Hardy's inequality (\ref{discreteHardy}).

\cite{MR4346672} prove (\ref{probabilityHardy}) with $\alpha=0$.
Their proof generalizes \cite{MR1574000}'s proof of (\ref{discreteHardy}), discretizes $F$ and applies a limiting procedure.
Theorem \ref{Hardytheorem} will be proved along the lines of (\ref{proofcontinuousHardy}) and to this end we need the following rearrangement lemma; cf. \cite{MR1574451}.

\begin{lemma}
Let $\chi$ be a nonnegative $p$-integrable function on the unit interval.
There exists a nonincreasing $p$-integrable function $\tilde \chi$ on the unit interval with the same $p$-norm as $\chi$ and with
\begin{equation}\label{chi}
\frac 1u \int_0^u {\tilde \chi}(v)\, dv \geq \frac 1u \int_0^u \chi(v)\, dv, \quad 0 \leq u \leq 1.
\end{equation}
\end{lemma}

\noindent 
{\bf Proof} \\
Let $U$ be uniformly distributed on the unit interval and define the distribution function $G(x) = P(\chi(U) \leq x), x \in {\mathbb R},$
and the nonincreasing function ${\tilde \chi}(u) = G^{-1}(1-u) = \inf \{ x \mid G(x) \geq 1-u \},\ 0 \leq u \leq 1.$
Note that ${\tilde \chi}(U)$ and $\chi(U)$ have the same distribution and hence the same $p$-norm.
Furthermore, by definition of $G$
\begin{equation*}
\frac 1u \int_0^u {\tilde \chi}(v)\, dv = \frac 1u \int_0^u G^{-1}(1-v)\, dv
\end{equation*}
is the mean of the $\chi$-values over a subset of $[0,1]$ of measure $u$ that contains the largest $\chi$-values.
This implies (\ref{chi}).
\hfill$\Box$ \\

Here is our proof of Theorem \ref{Hardytheorem}.

\noindent {\bf Proof} \\
Without loss of generality we assume that $\psi$ is nonnegative and $p$-integrable.
The left-continuous inverse distribution function $F^{-1}$ is defined as $F^{-1}(u) = \inf \{ x \mid F(x) \geq u \},\ 0 \leq u \leq 1.$
Since the left hand side of (\ref{probabilityHardy}) to the power $p$ may be rewritten as
\begin{equation}\label{rewrite}
\int_0^1 \left[ \frac { \int_{[F^{-1}(v) \leq F^{-1}(u)]} \psi(F^{-1}(v))\, dv }{F(F^{-1}(u))} \right]^p du
= \int_0^1 \left[ \frac { \int_0^{F(F^{-1}(u))} \psi(F^{-1}(v))\, dv }{F(F^{-1}(u))} \right]^p du,
\end{equation}
the lemma proves that we may assume that $\psi$ is nonincreasing.

Let $F$ be a distribution function with discontinuities, i.e., atoms, one of which is located at $a$ with mass $p_a$.
We remove this discontinuity by stretching the distribution function $F$ to
\begin{equation} \label{Ftilde}
{\tilde F}(x) = \left\{
    \begin{array}{lcl}
    F(x)         &                & x < a \\
    F(a-) + x -a & \mbox{\rm for} & a \leq x \leq a + p_a  \\
    F(x-p_a)     &                & a + p_a \leq x \\
    \end{array} \right.
\end{equation}
and adapt the function $\psi$ accordingly as follows
\begin{equation} \label{psitilde}
{\tilde \psi}(x) = \left\{
    \begin{array}{lcl}
    \psi(x)      &                & x < a \\
    \psi(a)      & \mbox{\rm for} & a \leq x \leq a + p_a  \\
    \psi(x-p_a)  &                & a + p_a \leq x.\\
    \end{array} \right.
\end{equation}
Note
\begin{equation}\label{equalitynorms}
\int_{\mathbb R} {\tilde \psi}^p(y)\, d{\tilde F}(y) = \int_{\mathbb R} \psi^p(y)\, dF(y).
\end{equation}
Since this also holds for $p=1$, the value of $\alpha$ in (\ref{root}) does not change
when replacing $\psi$ and $F$ by $\tilde \psi$ and $\tilde F$, respectively.
For $a \leq x \leq a + p_a$ the monotonicity of $\psi$ implies
\begin{eqnarray*}
\lefteqn{ \frac {\int_{(-\infty,x]} {\tilde \psi}(y)\, d{\tilde F}(y)}{ {\tilde F}(x)}
= \frac {\int_{(-\infty, a)} \psi(y)\, dF(y) + \psi(a) (x-a)}{F(a-) + x - a} } \\
&& = \frac {\int_{(-\infty, a)} [\psi(y) - \psi(a)] dF(y)}{F(a-) + x - a} + \psi(a) \\
&& \geq \frac {\int_{(-\infty, a]} [\psi(y) - \psi(a)] dF(y)}{F(a)} + \psi(a)
= \frac {\int_{(-\infty, a]} \psi(y)\, dF(y)}{F(a)}
\end{eqnarray*}
and for $a + p_a < x$
\begin{eqnarray*}
\lefteqn{ \frac {\int_{(-\infty,x]} {\tilde \psi}(y)\, d{\tilde F}(y)}{ {\tilde F}(x)} } \\
&& = \frac {\int_{(-\infty,a)} \psi(y)\, dF(y) + \psi(a)p_a + \int_{(a + p_a,x]} \psi(y-p_a)\,dF(y-p_a)}{F(x-p_a)} \\
&& = \frac {\int_{(-\infty, x-p_a]} \psi(y)\, dF(y)}{F(x-p_a)} 
\end{eqnarray*}
holds.

Together with the definitions of $\tilde F$ and $\tilde \psi$ these relations yield
\begin{eqnarray}\label{monotonicityinequality}
\lefteqn{ \int_{\mathbb R} \left[ \frac {\int_{(-\infty,x]} {\tilde \psi}(y)\, d{\tilde F}(y)}{ {\tilde F}(x)} \right]^p d{\tilde F}(x) \nonumber } \\
&& \geq \int_{(-\infty,a)} \left[ \frac {\int_{(-\infty,x]} \psi(y)\, dF(y)}{F(x)} \right]^p dF(x)
        + \int_a^{a+p_a} \left[ \frac {\int_{(-\infty, a]} \psi(y)\, dF(y)}{F(a)} \right]^p dx \nonumber \\
&&  \hspace{5em}  + \int_{(a+p_a,\infty)} \left[ \frac {\int_{(-\infty, x-p_a]} \psi(y)\, dF(y)}{F(x-p_a)} \right]^p dF(x-p_a) \\
&& = \int_{\mathbb R} \left[ \frac {\int_{(-\infty,x]} \psi(y)\, dF(y)}{F(x)} \right]^p dF(x). \nonumber
\end{eqnarray}
Since an arbitrary distribution function $F$ has at most countably many discontinuities, it might be that the stretch procedure from (\ref{Ftilde}) and (\ref{psitilde}) has to be repeated countably many times in order to obtain a continuous distribution function $\tilde F$ and an adapted function $\tilde \psi$ such that (\ref{equalitynorms}) holds and the left hand side of (\ref{monotonicityinequality}) equals at least its right hand side.
Consequently, we may and do assume that $F$ is continuous.

Next we prove
\begin{equation}\label{derivativeintegral}
\left[ \int_{(-\infty,x]} \psi \, dF \right]^p = p \int_{(-\infty,x]} \left[ \int_{(-\infty,y]} \psi \, dF \right]^{p -1} \psi(y) \, dF(y).
\end{equation}
Defining the distribution function
$F_x(y) =  \int_{(-\infty,y\wedge x]} \psi \, dF /  \int_{(-\infty,x]} \psi \, dF$ for $x$ with $F(x)>0$ 
we see that this equality is equivalent to
\begin{equation}\label{equalsone}
1 = p \int_{-\infty}^\infty F_x^{p -1}(y) \, dF_x(y) = p E\left( \left[F_x \left( F_x^{-1}(U) \right) \right]^{p -1} \right),
\end{equation}
where the random variable $U$ is uniformly distributed on the unit interval.
Since $F$ has no point masses, $F_x$ has none, i.e., $F_x$ is continuous.
Consequently $F_x(F_x^{-1}(u))=u,\, 0<u<1,$ holds and hence (\ref{equalsone}) and (\ref{derivativeintegral}). 
By (\ref{derivativeintegral}) and Tonelli's theorem we have
\begin{eqnarray}\label{boundLHS}
\lefteqn{ \int_{-\infty}^\infty \left( \frac {\int_{(-\infty, x]} \psi(y)\, dF(y)}{F(x)} \right)^p dF(x) \nonumber } \\
&& = p \int_{-\infty}^\infty (F(x))^{-p} \int_{(-\infty,x]} \left[ \int_{(-\infty,y]} \psi \, dF \right]^{p -1} \psi(y)\, dF(y)\, dF(x) \\
&& = p \int_{-\infty}^\infty \int_{[y,\infty)} (F(x))^{-p}\, dF(x) \left[ \int_{(-\infty,y]} \psi \, dF \right]^{p -1} \psi(y)\, dF(y). \nonumber
\end{eqnarray}
Since $F$ is continuous,
\begin{eqnarray*}
\lefteqn{ \int_{[y,\infty)} (F(x))^{-p}\, dF(x) = \int_{[y \leq F^{-1}(u)]} (F(F^{-1}(u)))^{-p}\, du \nonumber } \\
&& \hspace{10em} = \int_{F(y)}^1 u^{-p}\, du = \frac 1{p-1} \left[ (F(y))^{1-p} -1 \right]
\end{eqnarray*}
holds.
Together with (\ref{boundLHS}), H\"older's inequality and (\ref{derivativeintegral}) with $x=\infty$ we obtain
\begin{eqnarray*}
\lefteqn{ \int_{-\infty}^\infty \left( \frac {\int_{^(-\infty, x]} \psi(y)\, dF(y)}{F(x)} \right)^p dF(x) \nonumber } \\
&& = \frac p{p-1} \int_{-\infty}^\infty \left[ \frac 1{F(y)} \int_{(-\infty,y]} \psi \, dF \right]^{p -1} \psi(y)\, dF(y) \nonumber \\
&& \quad - \frac p{p-1} \int_{-\infty}^\infty \left[ \int_{(-\infty,y]} \psi \, dF \right]^{p -1} \psi(y) \, dF(y) \\
&&  \leq \frac p{p-1} \left[ \int_{-\infty}^\infty \left( \frac {\int_{^(-\infty, y]} \psi \, dF}{F(y)} \right)^p dF(y) \right]^{(p-1)/p}
  \left[ \int_{-\infty}^\infty \psi^p(y)\, dF(y) \right]^{1/p} \nonumber \\
&& \quad - \frac 1{p-1} \left[ \int_{-\infty}^\infty \psi(y) \, dF(y) \right]^p.  \nonumber
\end{eqnarray*}
Introducing shorthand notation we write this inequality as
\begin{equation*}
z^p \leq \frac p{p-1} z^{p-1} \nu - \frac 1{p-1} \mu^p.
\end{equation*}
This is equivalent to
\begin{equation}\label{chi0}
\chi(z) = p \nu z^{p-1} -(p-1)z^p - \mu^p \geq 0.
\end{equation}
The function $\chi$ is increasing -- decreasing on the positive half line with a nonnegative maximum at $\nu$.
Furthermore
\begin{equation}\label{chi1}
\chi \left( p\nu / (p-1+\alpha) \right) = 0
\end{equation}
holds in view of (\ref{root}).
Since the left hand side of (\ref{root}) is convex in $\alpha$ on $[0,\infty)$ and is nonnegative at $\alpha=0$ and nonpositive at $\alpha=1$,
the root $\alpha$ from (\ref{root}) satisfies $\alpha \in [0,1]$ and hence $\nu \leq  p\nu / (p-1+\alpha)$.
Together with (\ref{chi0}) and (\ref{chi1}) this proves \\ $z \leq  p\nu / (p-1+\alpha)$, i.e., (\ref{probabilityHardy}).
\hfill$\Box$ \\

\begin{remark}
Since the function
\begin{equation*}
u \mapsto \frac 1u \int_0^u \psi(F^{-1}(v)) dv, \quad 0 \leq u \leq 1,
\end{equation*}
is nonincreasing for nonnegative, nonincreasing $\psi$, we have
\begin{equation}\label{inequality}
\int_0^1 \left[ \frac { \int_0^{F(F^{-1}(u))} \psi(F^{-1}(v))\, dv }{F(F^{-1}(u))} \right]^p du
\leq \int_0^1 \left[ \frac 1u \int_0^u \psi(F^{-1}(v))\, dv \right]^p du
\end{equation}
in view of $F(F^{-1}(u)) \geq u$.
We define $\psi(F^{-1}(v))=0$ for $v>1$.
Applying (\ref{rewrite}), (\ref{inequality}) and (\ref{continuousHardy}) for nonnegative, nonincreasing $\psi$ we arrive at
\begin{eqnarray*}
\lefteqn{ E\left( \left[ \frac{E\left( \psi(Y) {\bf 1}_{[Y \leq X]} \mid X \right)}{F(X)} \right]^p \right)
\leq \int_0^1 \left[ \frac 1u \int_0^u \psi(F^{-1}(v))\, dv \right]^p du \nonumber } \\
&& \leq \int_0^\infty \left[ \frac 1u \int_0^u \psi(F^{-1}(v))\, dv \right]^p du
\leq \left(\frac p{p - 1} \right)^p \int_0^\infty \left( \psi(F^{-1}(v)) \right)^p dv \\
&& = \left(\frac p{p - 1} \right)^p \int_0^1 \left( \psi(F^{-1}(v)) \right)^p dv
= \left(\frac p{p - 1} \right)^p  E\left( \psi^p(Y) \right). \nonumber
\end{eqnarray*}
This proves (\ref{probabilityHardy}) with $\alpha=0$, which is the probability theoretic version of Hardy's inequality given in Theorem 2.1 of \cite{MR4346672}.
Arguably this proof is simpler and more elegant than the one given in \cite{MR4346672}.
\end{remark}

\section{A probability theoretic Hardy inequality for $0<p<1$}

Theorem 337 on page 251 of \cite{MR0046395} states that for $0<p<1$ and $\psi$ a nonnegative measurable function on $(0, \infty)$
\begin{eqnarray}\label{Hardy<1}
\left\{\int_0^\infty \left ( \frac{1}{x} \int_x^\infty \psi (y) dy \right )^p \, dx \right\}^{1/p}
\geq \frac p{1-p} \left\{ \int_0^{\infty} \psi^p(y)\, dy \right\}^{1/p}
\end{eqnarray}
holds.
A smoothed version of their proof similar to the one in (\ref{proofcontinuousHardy}) can be given, but with the inverse H\"older inequality this time;
\begin{eqnarray}\label{proofcontinuousHardy<1}
\lefteqn{ \int_0^\infty \left ( \frac{1}{x} \int_x^\infty \psi \right )^p \, dx
= \int_0^\infty x^{-p} \int_x^\infty p \left(\int_y^\infty \psi \right)^{p-1} \psi(y)\, dy \, dx \nonumber }\\
&& = p \int_0^\infty \int_0^y x^{-p} \, dx \left( \int_y^\infty \psi \right)^{p-1} \psi(y)\, dy \\
&&  = \frac p{1-p} \int_0^\infty \left( \frac 1y \int_y^\infty \psi \right)^{p-1} \psi(y) \, dy  \nonumber \\
&& \geq \frac p{1-p} \left[\int_0^\infty \left( \frac 1y \int_y^\infty \psi \right)^p dy \right]^{(p-1)/p}
\left[ \int_0^\infty \psi^p \right]^{1/p}. \nonumber
\end{eqnarray}
Our probability theoretic generalization of this inequality reads as follows.
\begin{theorem}\label{Hardytheorem<1}
Let $X$ and $Y$ be independent random variables with distribution function $F$ on $({\mathbb R}, {\cal B})$,
and let $\psi $ be a nonnegative measurable function on $({\mathbb R}, {\cal B})$.
For $0<p<1$
\begin{equation}\label{probabilityHardy<1}
\left\{ E\left( \left( \frac{E\left( \psi(Y) {\bf 1}_{[Y \geq X]} \mid X \right)}{F(X-)} \right)^p \right) \right\}^{1/p}
\geq \frac p{1-p} \left\{ E\left( \psi^p(Y) \right) \right\}^{1/p}
\end{equation}
holds.
\end{theorem}
\noindent{\bf Proof}
The proof is analogous to the one for Theorem \ref{Hardytheorem}, but with essential modifications.
Let $F$ be a distribution function with discontinuities, one of which is located at $a$ with jump size $p_a$.
We remove this discontinuity by stretching the distribution function $F$ to
\begin{equation} \label{Fbar}
{\bar F}(x) = \left\{
    \begin{array}{lcl}
    F(x+p_a)      &                & x < a - p_a \\
    F(a) + x -a   & \mbox{\rm for} & a - p_a \leq x \leq a \\
    F(x)          &                & a < x \\
    \end{array} \right.
\end{equation}
and adapt the function $\psi$ accordingly as follows
\begin{equation} \label{psibar}
{\bar \psi}(x) = \left\{
    \begin{array}{lcl}
    \psi(x + p_a) &                & x < a - p_a \\
    \psi(a)       & \mbox{\rm for} & a - p_a \leq x \leq a \\
    \psi(x)       &                & a < x.\\
    \end{array} \right.
\end{equation}
Note
\begin{equation}\label{equalitynormsbar}
\int_{\mathbb R} {\bar \psi}^p(y)\, d{\bar F}(y) = \int_{\mathbb R} \psi^p(y)\, dF(y).
\end{equation}
For $x < a - p_a$ we have
\begin{eqnarray*}
\lefteqn{ \frac {\int_{[x,\infty)} {\bar \psi}(y)\, d{\bar F}(y)}{ {\bar F}(x-)} } \\
&& = \frac {\int_{[x,a-p_a)} \psi(y+p_a)\, dF(y+p_a) + \psi(a)p_a + \int_{(a,\infty)} \psi(y)\, dF(y)}{F(x+p_a-)} \\
&& = \frac {\int_{[x+p_a,\infty)} \psi(y)\, dF(y)}{F(x+p_a-)},
\end{eqnarray*}
for $a - p_a \leq x \leq a$ the nonnegativity of $\psi$ implies
\begin{eqnarray*}
\frac {\int_{[x,\infty)} {\bar \psi}(y)\, d{\bar F}(y)}{ {\bar F}(x-)} 
= \frac {\psi(a)(a-x) + \int_{(a,\infty)} \psi(y)\, dF(y)}{ F(a) + x - a} 
\leq \frac {\int_{[a,\infty)} \psi(y)\, dF(y)}{F(a-)}
\end{eqnarray*}
and for $a<x$ we have
\begin{eqnarray*}
\int_{[x,\infty)} {\bar \psi}(y)\, d{\bar F}(y) / {\bar F}(x-)
= \int_{[x,\infty)} \psi(y)\, dF(y) / F(x-).
\end{eqnarray*}

Together with the definitions of $\bar F$ and $\bar \psi$ these relations yield
\begin{eqnarray}\label{nonnegativityinequality}
\lefteqn{ \int_{\mathbb R} \left[ \frac {\int_{[x,\infty)} {\bar \psi}(y)\, d{\bar F}(y)}{ {\bar F}(x-)} \right]^p d{\bar F}(x) \nonumber } \\
&& \leq \int_{(-\infty,a-p_a)} \left[ \frac {\int_{[x+p_a,\infty)} \psi(y)\, dF(y)}{F(x+p_a-)} \right]^p dF(x+p_a) \\
&& \hspace{2em} + \int_{a-p_a}^a \left[ \frac {\int_{[a,\infty)} \psi(y)\, dF(y)}{F(a-)} \right]^p dx
    + \int_{(a,\infty)} \left[ \frac {\int_{[x,\infty)} \psi(y)\, dF(y)}{F(x-)} \right]^p dF(x) \nonumber \\
&& = \int_{\mathbb R} \left[ \frac {\int_{[x,\infty)} \psi(y)\, dF(y)}{F(x-)} \right]^p dF(x). \nonumber
\end{eqnarray}
Since an arbitrary distribution function $F$ has at most countably many discontinuities, it might be that the stretch procedure from (\ref{Fbar}) and (\ref{psibar}) has to be repeated countably many times in order to obtain a continuous distribution function $\bar F$ and an adapted function $\bar \psi$ such that (\ref{equalitynormsbar}) holds and the left hand side of (\ref{nonnegativityinequality}) equals at most its right hand side.
Consequently, we may and do assume that $F$ is continuous.

For $x$ with $\int_{[x,\infty)} \psi \, dF >0$ we define
\begin{equation*}
G_x(y) = \int_{[x, y \vee x]} \psi \, dF \, / \int_{[x,\infty)} \psi \, dF,\quad y \in {\mathbb R}.
\end{equation*}
Since $F$ is continuous, $G_x$ is and we have
\begin{eqnarray*}
\lefteqn{ p \int_{\mathbb R} \left( 1-G_x(y) \right)^{p-1} dG_x(y) = p \int_0^1 \left( 1 -G_x(G_x^{-1}(u)) \right)^{p-1} du } \\
&& \hspace{10em} = p \int_0^1 (1-u)^{p-1} du = 1
\end{eqnarray*}
and hence
\begin{equation}\label{derivativeintegralbar}
\left( \int_{[x,\infty)} \psi \, dF \right)^p = p \int_{[x,\infty)} \left( \int_{[y,\infty)} \psi \, dF \right)^{p-1} \psi(y) \, dF(y).
\end{equation}
By (\ref{derivativeintegralbar}), Tonelli's theorem, the continuity of $F$ and the inverse H\"older inequality we obtain
\begin{eqnarray*}
\lefteqn{ \int_{-\infty}^\infty \left( \frac {\int_{[x,\infty)} \psi(y)\, dF(y)}{F(x-)} \right)^p dF(x) \nonumber } \\
&& = p \int_{-\infty}^\infty (F(x))^{-p} \int_x^\infty \left[ \int_y^\infty \psi \, dF \right]^{p -1} \psi(y)\, dF(y)\, dF(x) \nonumber \\
&& = p \int_{-\infty}^\infty \int_{-\infty}^y (F(x))^{-p} dF(x) \left[ \int_y^\infty \psi \, dF \right]^{p -1} \psi(y)\, dF(y) \\
&& = \frac p{1-p} \int_{-\infty}^\infty \left( \frac {\int_y^\infty \psi \, dF}{F(y)} \right)^{p-1} \psi(y)\, dF(y) \nonumber \\
&& \geq \frac p{1-p} \left[ \int_{-\infty}^\infty \left( \frac {\int_y^\infty \psi \, dF}{F(y)} \right)^p \, dF(y) \right]^{(p-1)/p}
\left[ \int_{-\infty}^\infty \psi^p(y) \, dF(y) \right]^{1/p},  \nonumber
\end{eqnarray*}
which implies (\ref{probabilityHardy<1}).
\hfill$\Box$ \\

Note that taking $F$ uniform on $(0,K)$, multiplying (\ref{probabilityHardy<1}) by $K$
and taking limits as $K \to \infty$ we arrive at Hardy's inequality (\ref{Hardy<1}).

A discrete version of Theorem 337 from \cite{MR0046395} and its proof are given in Theorem 338 on page 252 of \cite{MR0046395}.
This Theorem 338 states that for $0<p<1$ and $a_i \geq 0, \ i=1, 2, \dots,$
\begin{equation}\label{discreteHardy<1}
\left( 1 + \frac 1{1-p} \right) \left( \sum_{i=1}^\infty a_i \right)^p + \sum_{j=2}^\infty \left( \frac 1j \sum_{h=j}^\infty a_h \right)^p 
\geq \left( \frac p{1-p} \right)^p \sum_{i=1}^\infty a_i^p.
\end{equation}
This inequality follows from our Theorem \ref{Hardytheorem<1} by the choices
\begin{equation}\label{choices}
F(x) = \frac xK {\bf 1}_{[0<x<1]} + \sum_{i=1}^K \frac 1K {\bf 1}_{[i \leq x]}, \ \psi(x) = \sum_{i=1}^\infty a_i {\bf 1}_{[x=i+1]}, \ x \in {\mathbb R},
\end{equation}
where $K$ is a natural number.
Indeed, we have
\begin{equation}\label{limit1}
\lim_{K \to \infty} K \  E\left( \psi^p(Y) \right) = \lim_{K \to \infty} \sum_{i=1}^{K-1} a_i^p = \sum_{i=1}^\infty a_i^p
\end{equation}
and
\begin{eqnarray}\label{limit2}
\lefteqn{\lim_{K \to \infty} K \  E\left( \left( \frac{E\left( \psi(Y) {\bf 1}_{[Y \geq X]} \mid X \right)}{F(X-)} \right)^p \right) \nonumber } \\
&& = \lim_{K \to \infty} \int_0^1 \left( \frac {\sum_{i=1}^{K-1} a_i }x \right)^p dx 
      + \sum_{j=2}^K \left( \frac 1{j-1} \sum_{h=j-1}^{K-1} a_h \right)^p \\
&& = \frac 1{1-p}  \left( \sum_{i=1}^\infty a_i \right)^p + \sum_{h=1}^\infty \left( \frac 1h \sum_{j=h}^\infty a_j \right)^p, \nonumber
\end{eqnarray}
which equals the left hand side of (\ref{discreteHardy<1}).

\begin{remark}
For $p=1$ no general inequality exists like in Theorems \ref{Hardytheorem} and \ref{Hardytheorem<1}.
With $\psi$ nonnegative nondecreasing we have
\begin{equation*}
E\left( \frac{E\left( \psi(Y) {\bf 1}_{[Y \leq X]} \mid X \right)}{F(X)} \right)
\leq E\left( \frac{\psi(X) {\bf 1}_{[Y \leq X]}}{F(X)} \right) = E(\psi(X))
\end{equation*}
and with $\psi$ nonnegative nonincreasing
\begin{equation*}
E\left( \frac{E\left( \psi(Y) {\bf 1}_{[Y \leq X]} \mid X \right)}{F(X)} \right)
\geq E\left( \frac{\psi(X) {\bf 1}_{[Y \leq X]}}{F(X)} \right) = E(\psi(X)).
\end{equation*}
Similarly we get $E(\psi(X)(1-F(X))/F(X-))$ as both an upper and a lower bound to the left hand side of (\ref{probabilityHardy<1}) with $p=1$.
\end{remark}

\section{The probability theoretic Copson inequality}

Dual to Hardy's inequality is Copson's inequality.
It was presented in \cite{MR1574056}; see also Section 5 of \cite{MR4346672}.
Inspired by (\ref{proofcontinuousHardy}) we present a short proof of the probability theoretic Copson inequality as given in \cite{MR4346672}.
\begin{theorem}\label{Copson}
Let $X$ and $Y$ be independent random variables with distribution function $F$ on $({\mathbb R}, {\cal B})$,
and let $\psi $ be a measurable function on $({\mathbb R}, {\cal B})$.
With $p \geq 1$
\begin{equation}\label{probabilityCopson}
\left\{ E\left( \left| E\left( \frac {\psi(Y)}{F(Y)} {\bf 1}_{[Y \geq X]} \mid X \right) \right|^p \right) \right\}^{1/p}
\leq p \, \left\{ E\left( |\psi(Y)|^p \right) \right\}^{1/p}
\end{equation}
holds.
\end{theorem}
\noindent {\bf Proof} \\
For $p=1$ we have
\begin{equation*}
E\left( \left| E\left( \frac {\psi(Y)}{F(Y)} {\bf 1}_{[Y \geq X]} \mid X \right) \right| \right) 
\leq E\left( \frac {|\psi(Y)|}{F(Y)} {\bf 1}_{[Y \geq X]} \right) = E( |\psi(Y)| ).
\end{equation*}
Without loss of generality we assume $p>1$ and that $\psi$ is nonnegative and $p$-integrable.
For $x$ with $F(x)>0$ we have $\int_{[x,\infty)}\psi/F \, dF \leq E(\psi(X))/F(x) < \infty$.
Consequently, for such $x$ that satisfy $\int_{[x,\infty)} \psi / F\, dF >0$ as well
\begin{equation*}
G_x(y) = \int_{[x,y]} \frac {\psi(z)}{F(z)} dF(z) \left[\int_{[x,\infty)} \frac {\psi(z)}{F(z)} dF(z) \right]^{-1}, \quad y \in {\mathbb R},
\end{equation*}
is a well defined distribution function.
For any distribution function $G$ and corresponding left-continuous inverse distribution function $G^{-1}$
the inequalities $G(G^{-1}(u)) \geq u$ and $G(G^{-1}(u)-) \leq u$ hold.
Consequently with $U$ uniformly distributed we have
\begin{eqnarray*}
\lefteqn{ p \, E_{G_x} \left( \left[ 1 - G_x(Y-) \right]^{p-1} \right)
= p \, E \left( \left[ 1 - G_x \left( G_x^{-1}(U) - \right) \right]^{p-1} \right) } \\
&& \hspace{10em} \geq p \, E \left( [1-U]^{p-1} \right) = 1.
\end{eqnarray*}
Combining this inequality with Tonelli's theorem and H\" older's inequality we see that the left hand side of (\ref{probabilityCopson}) equals and satisfies
\begin{eqnarray*}
\lefteqn{ \int_{-\infty}^\infty \left[ \int_{[x, \infty)} \frac \psi F \, dF \right]^p dF(x) } \\
&& \leq p \int_{-\infty}^\infty \int_{[x, \infty)} \left[ \int_{[y, \infty)} \frac \psi F \, dF \right]^{p-1} \frac {\psi(y)}{F(y)} \, dF(y) \ dF(x) \\
&& =  p \int_{-\infty}^\infty \int_{(-\infty,y]} dF(x) \left[ \int_{[y, \infty)} \frac \psi F \, dF \right]^{p-1} \frac {\psi(y)}{F(y)} \, dF(y) \\
&& =  p \int_{-\infty}^\infty \left[ \int_{[y, \infty)} \frac \psi F \, dF \right]^{p-1} \psi(y) \, dF(y) \\
&& \leq p \left\{ \int_{-\infty}^\infty \left[ \int_{[y, \infty)} \frac \psi F \, dF \right]^p dF(y) \right\}^{(p-1)/p}
\left\{ \int_{-\infty}^\infty \psi^p(y) \, dF(y) \right\}^{1/p},
\end{eqnarray*}
which proves (\ref{probabilityCopson}).
\hfill$\Box$ 

Slight modifications in this proof yield a probability theoretic generalization of Copson's inequality for $0<p<1$; cf. Theorem 2.3 with $c=\kappa=p$ in \cite{MR1574443}.
\begin{theorem}\label{Copson<1}
Let $X$ and $Y$ be independent random variables with distribution function $F$ on $({\mathbb R}, {\cal B})$,
and let $\psi $ be a nonnegative measurable function on $({\mathbb R}, {\cal B})$.
With $0<p<1$
\begin{equation}\label{probabilityCopson<1}
\left\{ E\left( \left( E\left( \frac {\psi(Y)}{F(Y)} {\bf 1}_{[Y \geq X]} \mid X \right) \right)^p \right) \right\}^{1/p}
\geq p \, \left\{ E\left( \psi^p(Y) \right) \right\}^{1/p}
\end{equation}
holds.
\end{theorem}
\noindent {\bf Proof} \\
We may and do assume that the left hand side of (\ref{probabilityCopson<1}) is finite and hence that $\int_{[x,\infty)}\psi/F \, dF <\infty$ holds for $F$-almost all $x$.
Consequently, for such $x$ with $F(x)>0$ and $\int_{[x,\infty)}\psi/F \, dF > 0$ the distribution function $G_x$ from the preceding proof is well defined.
Noting $0<p<1$ and replacing H\"older's inequality by its reversed version we see that the preceding proof is valid with all inequality signs reversed.
\hfill$\Box$ \\
\vspace{1cm} \\
\noindent
{\bf Acknowledgements}
The author would like to thank Jon Wellner for suggesting the topic of extending Hardy's inequality in a probability theoretic way, for his collaboration on this topic resulting in \cite{MR4346672}, and for suggesting to write this note.
Furthermore I owe thanks to the late Jaap Fabius who once said to me: "Partial integration is just Fubini."

\end{document}